# Jack characters and enumeration of maps

Agnieszka Czyżewska-Jankowska[1] and Piotr Śniady[2]

[1] *Piławska 20/16, 50-538 Wrocław, Poland* [2] *Institute of Mathematics, Polish Academy of Sciences, ul. Śniadeckich 8, 00-956 Warszawa, Poland*

**Abstract.** *Jack characters* provide dual information about Jack symmetric functions. We give explicit formulas for the top-degree part of these Jack characters in terms of bicolored oriented maps with an arbitrary face structure.

## 1 Jack characters are interesting because...

This short note presents ideas from preprints [3, 11] which will be published elsewhere.

For an integer partition $\pi$, *Jack character* $\mathrm{Ch}_\pi$ is a certain function on the set $\mathbb{Y}$ of Young diagrams, valued in Laurent polynomials $\mathbb{Q}[A, A^{-1}]$. For example,

$$\mathrm{Ch}_3(\lambda) = \sum_{\square_1 \in \lambda} \left( 3(c_1 + \gamma)(c_1 + 2\gamma) + \frac{3}{2} \right) + \sum_{\square_1, \square_2 \in \lambda} \left( -\frac{3}{2} \right), \quad (1.1)$$

where the sums run over the boxes of the Young diagram $\lambda$ and we use the notation that

$$c_i = c_i(\square_i) = A x_i - A^{-1} y_i \quad (1.2)$$

denotes the *A-deformed content* of a box $\square_i = (x_i, y_i) \in \mathbb{N}^2$ which is in $x_i$-th column and $y_i$-th row (cf. Figure 2b), and

$$\gamma = -A + A^{-1}. \quad (1.3)$$

*There are five good reasons for studying Jack characters. We shall review them in the following.*

### 1.1 ...they are related to Jack polynomials.

*Jack characters provide dual information about Jack polynomials* which are a 'simple' version of Macdonald polynomials. Despite this 'simplicity', their structure remains elusive [8]. *Therefore a better understanding of Jack characters might shed some light on Jack polynomials.*

To be more specific, we expand the Jack polynomial $J_\pi^{(\alpha)}$ which corresponds to the deformation parameter $\alpha := A^2$ in the basis of power-sum symmetric functions:

$$J_\lambda^{(\alpha)} = \sum_\pi \theta_\pi^{(\alpha)}(\lambda) \, p_\pi.$$

Jack character $\mathrm{Ch}_\pi(\lambda)$ is equal to suitably normalized coefficient $\theta_\pi^{(\alpha)}(\lambda)$, cf. Lassalle [9] (with another normalization) and Dołęga and Féray [4] (whose normalization we use).

## 1.2  . . . they are a deformation of the symmetric group characters. . .

. . . which depends on an additional parameter $A$. In enumerative combinatorics it is quite common that such deformation parameters might shed some light into the structure of the original, non-deformed object. Therefore *a better understanding of Jack characters might be beneficial also for the investigation of the irreducible characters of the symmetric groups.*

To be more concrete: for the specific choice $A=1$ each Jack character coincides with the (suitably normalized) usual character of the symmetric group:

$$\mathrm{Ch}_\pi(\lambda)\Big|_{A:=1} = \underbrace{|\lambda|\cdot(|\lambda|-1)\cdots(|\lambda|-|\pi|+1)}_{|\pi|\text{ factors}} \frac{\mathrm{Tr}\,\rho^\lambda(\pi, 1^{|\lambda|-|\pi|})}{\mathrm{Tr}\,\rho^\lambda(1^{|\lambda|})}. \qquad (1.4)$$

## 1.3  . . . they have interesting structure constants.

We expand a product of two Jack characters in the basis provided by Jack characters, e.g.

$$\mathrm{Ch}_3\,\mathrm{Ch}_3 = (6\delta^2+3)\,\mathrm{Ch}_3 + 9\delta\,\mathrm{Ch}_{2,1} + 18\delta\,\mathrm{Ch}_4 + 3\,\mathrm{Ch}_{1,1,1} + 9\,\mathrm{Ch}_{3,1} + 9\,\mathrm{Ch}_{2,2} + 9\,\mathrm{Ch}_5 + \mathrm{Ch}_{3,3}.$$

The coefficients in such expansions are conjecturally polynomials in the variable $\delta = -\gamma$ (cf. (1.3)) with *nonnegative integer coefficients*, the combinatorial meaning of which remains elusive but indicates some natural deformation of the symmetric group algebra [10].

## 1.4  . . . there are interesting formulas for them.

Specifically: consider $A, p_1, \ldots, p_\ell, q_1, \ldots, q_\ell \in \mathbb{R}$ with the property that the shape depicted on Figure 1 defines a Young diagram. It turns out that the value of $\mathrm{Ch}_\pi$ on such a diagram is a polynomial (called *Stanley character polynomial*) in the variables $\gamma, p_1, \ldots, p_\ell, q_1, \ldots, q_\ell$, cf. (1.3). For example, in the case $\ell=2$ of two rectangles:

$$-\mathrm{Ch}_3 = \qquad\qquad\qquad\qquad\qquad\qquad\qquad\qquad\qquad\qquad\qquad\qquad (1.5)$$
$$\left.\begin{array}{l} p_1^3 q_1 + 3p_1^2 q_1^2 + p_1 q_1^3 + 3p_1^2 p_2 q_2 + 3p_1 p_2^2 q_2 + p_2^3 q_2 + 3p_1 p_2 q_1 q_2 + 3p_1 p_2 q_2^2 + 3p_2^2 q_2^2 + \\ p_2 q_2^3 + 3p_1^2 q_1 \gamma + 3p_1 q_1^2 \gamma + 6p_1 p_2 q_2 \gamma + 3p_2^2 q_2 \gamma + 3p_2 q_2^2 \gamma + 2p_1 q_1 \gamma^2 + 2p_2 q_2 \gamma^2 \end{array}\right\} =: \mathrm{Ch}_3^{\mathrm{top}}$$
$$+ p_1 q_1 + p_2 q_2.$$

**Coefficients of such polynomials for $\mathrm{Ch}_n$ seem to be** (up to a global sign change) **nonnegative integers. We will show in this note that this is indeed the case for $\mathrm{Ch}_n^{\mathrm{top}}$, the top-degree homogeneous part, indicated above by the curly brace.**

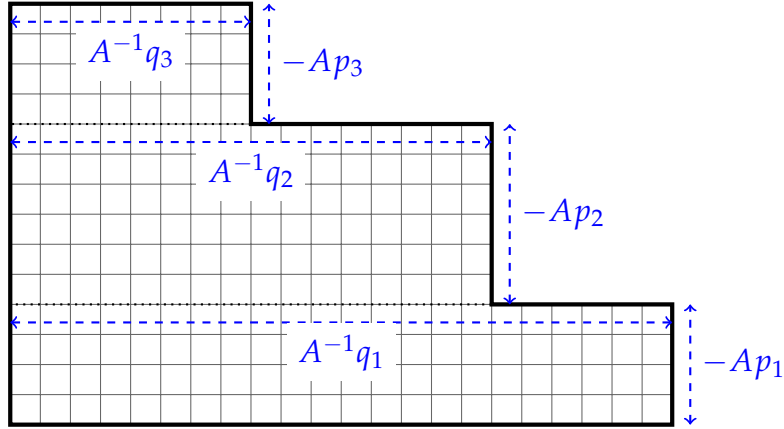

Figure 1: Multirectangular Young diagram $(-A\mathbf{p}) \times (A^{-1}\mathbf{q})$.

## 1.5 ...they can be defined in a convenient, an abstract way...

...which does not refer to the notion of Jack polynomials. Bear in mind example (1.1).

*Definition* 1.1. We say that $F\colon \mathbb{Y} \to \mathbb{Q}\left[A, A^{-1}\right]$ is a *polynomial function of degree at most* $d \in \{0, 1, \dots\}$ if there exist polynomials $p_0, p_1, \dots, p_{\lfloor \frac{d}{2} \rfloor}$ such that:

(J1) for each $0 \leq k \leq \lfloor \frac{d}{2} \rfloor$ we have that $p_k \in \mathbb{Q}[\gamma, c_1, \dots, c_k]$ is of degree at most $d - 2k$ and $p_k$ regarded as a polynomial in $c_1, \dots, c_k$ with coefficients in $\mathbb{Q}[\gamma]$ is symmetric;

(J2) for each Young diagram $\lambda$,

$$F(\lambda) = \sum_{0 \leq k \leq \lfloor \frac{d}{2} \rfloor} \sum_{\square_1, \dots, \square_k \in \lambda} p_k(\gamma, c_1, \dots, c_k) \in \mathbb{Q}\left[A, A^{-1}\right], \tag{1.6}$$

where the conventions (1.2), (1.3) are used.

We concentrate on Jack characters related to partitions $\pi = (n)$ with a single part.

*Definition* 1.2. $\mathrm{Ch}_n$ is the unique polynomial function of degree $d := n + 1$ such that:

(J3)
$$\mathrm{Ch}_n(\lambda) = 0 \qquad \text{for each } \lambda \in \mathbb{Y} \text{ such that } |\lambda| < n;$$

(J4) the top-degree coefficient of $p_1$ (from Definition 1.1) with respect to $c_1$ is given by

$$\left[c_1^{n-1}\right] p_1 = n.$$

The meaning of the above conditions can be explained heuristically, at least in the special case $A = 1$, $\gamma = 0$ (cf. Section 1.2) of the characters of the symmetric groups. (J2) and existence of the content polynomials $p_0, p_1, \ldots$ was proved by Corteel, Goupil and Schaeffer [2] and is related to *Jucys–Murphy elements*. (J1) and the degree bounds on the content polynomials on one hand, and (J4) on the other hand, reflect the asymptotics of the characters of the symmetric groups $\text{Ch}_\pi(\lambda)|_{A=1}$ in the scaling of *balanced Young diagrams* [1] (i.e. $\lambda \to \infty$ and its number of rows and columns grows like $O(\sqrt{|\lambda|})$) on one hand, and in the *Thoma scaling* [12] (in which the number of rows and columns of $\lambda$ grows like $\Theta(|\lambda|)$) on the other hand. (J3) reflects the normalization factor in (1.4) which vanishes on small Young diagrams.

With this in mind, the Jack character $\text{Ch}_\pi$ is indeed a natural generalization of the usual characters of the symmetric groups, with the only difference that the notion of content of a box is replaced by its deformation, the *A-content* (1.2).

## 2 Characters $\text{Ch}_n$ and maps

### 2.1 How to prove a formula for $\text{Ch}_n$?

Definition 1.2 opens the opportunity of proving a formula for Jack character in the following two easy steps: (1) guess a closed formula for $\text{Ch}_n$; then (2) verify that the guessed formula fulfills the defining properties of Jack character. Regretfully, already the first easy step is a challenge. We shall review the attempts to overcome it.

### 2.2 Maps

A *map* is a graph $G$ drawn on a surface $S$. Each map which we shall consider today is: *bicolored*, i.e. it is bipartite with a specific choice of the decoposition $\mathcal{V} = \mathcal{V}_\circ \sqcup \mathcal{V}_\bullet$ of the set of vertices into white and black vertices; *connected*; *rooted*, i.e. one of the edges is decorated; *unlabeled*, i.e. the remaining edges do not carry any additional decorations. For an example see Figure 2a with all labels of the vertices removed.

We shall consider two classes of such maps: a map might be *oriented* if the surface $S$ is orientable and it comes with a specified choice of the orientation; or it might be *not oriented* if we make no assumptions about orientability of $S$.

### 2.3 Embeddings of graphs

A pair of functions $(f_1, f_2)$ is called an *embedding* of a map $M$ to a Young diagram $\lambda$ if

$$f_1 \colon \mathcal{V}_\circ \to \text{(the set of columns of } \lambda\text{)}, \qquad f_2 \colon \mathcal{V}_\bullet \to \text{(the set of rows of } \lambda\text{)}$$

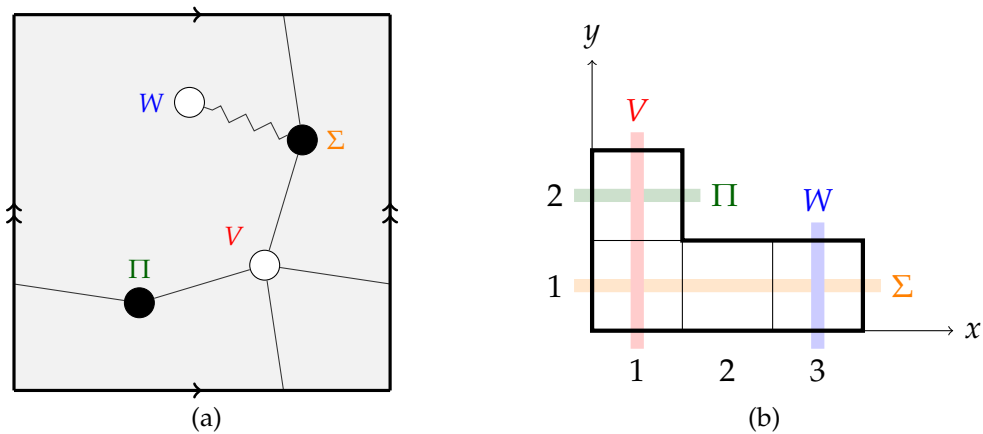

Figure 2: (a) Map on the torus and (b) an example of its embedding $f_1(V) = 1$, $f_2(W) = 3$, $f_2(\Sigma) = 1$, $f_2(\Pi) = 2$.

are such that for all $w \in \mathcal{V}_\circ$, $b \in \mathcal{V}_\bullet$ connected by an edge, $(f_1(w), f_2(b)) \in \lambda$, i.e. the intersection of column $f_1(w)$ and the row $f_2(b)$ belongs to $\lambda$, cf. Figure 2b. We set

$$\mathfrak{N}_G(\lambda) := A^{|\mathcal{V}_\circ(G)|} (-A)^{-|\mathcal{V}_\bullet(G)|} \times (\text{the number of embeddings of } G \text{ to } \lambda) \in \mathbb{Q}\left[A, A^{-1}\right].$$

## 2.4 Special values of $A$

For $A = 1$ a closed formula for the characters of the symmetric groups is available [7]:

$$\mathrm{Ch}_n \Big|_{A:=1} = (-1) \sum_M \mathfrak{N}_M, \tag{2.1}$$

where the sum runs over **oriented** maps $M$ with $n$ edges **and one face**.

An analogue of (2.1) holds true also for $A \in \left\{\sqrt{2}, \frac{1}{\sqrt{2}}\right\}$; the only difference (apart from some simple numerical factor) is that the sum on the right-hand side runs over **non-oriented** maps $M$ with $n$ edges **and one face**.

## 2.5 Great expectations and their depressive end

The above examples indicate existence of some hypothetical formula for $\mathrm{Ch}_n$ in the generic case which would be analogous to (2.1), with the sum running over **non-oriented** maps $M$ with $n$ edges **and one face**; each summand should be multiplied by some weight depending on $\gamma$ which would measure the 'orientability' of the map $M$.

A concrete form of this hypothetical 'measure of non-orientability' has been proposed in [5]. One should erase the edges from the *ribbon graph* of a map (cf. Figure 3) one after another in some random linear order $\prec$. To an edge $e$ which is to be removed one associates a factor which depends on the way the edge $e$ is attached to the remaining part of the ribbon graph. Roughly speaking, one regards whether the edge $e$ is a part of a Möbius strip, and whether the two sides of $e$ belong to the same face or not. The 'measure of non-orientability' was defined in [5] as the mean value (over the choice of $\prec$) of the product of such factors.

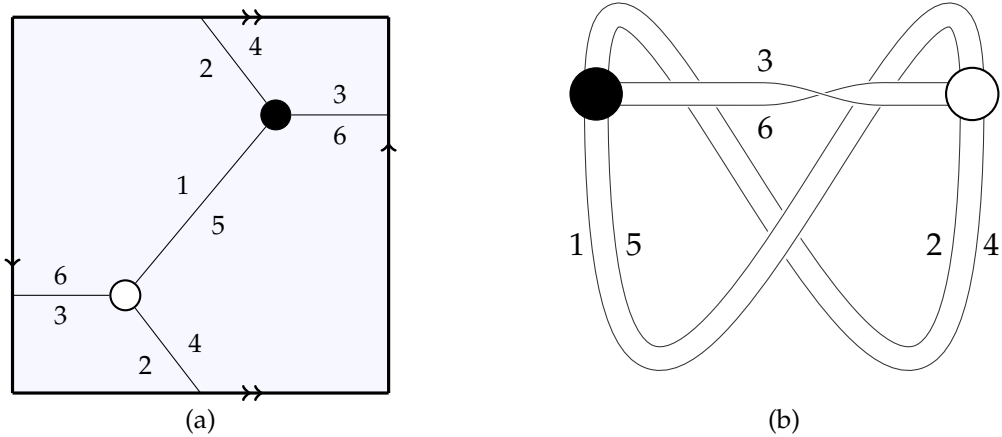

Figure 3: (a) Example of a map drawn on the Klein bottle and (b) the corresponding *ribbon graph*.

Regretfully, *the conjectural formula from [5] turned out to be incorrect*, which might seem a killing blow in the step (1) of our plan from Section 2.1 and the end of the story.

## 2.6   Hero of the day: top degree part of $\mathrm{Ch}_n$

Nevertheless, the conjectural formula from [5] turned out to predict *some* properties of Jack characters *surprisingly well*. For example, computer experiments indicated that it gives the correct value for the *most* of the coefficients of Stanley polynomials, cf. (1.5). Maybe we can take advantage of this observation and achieve a more modest but more realistic goal: *prove a closed formula for the homogeneous top-degree part $\mathrm{Ch}_n^{\mathrm{top}}$ of the Stanley polynomial for $\mathrm{Ch}_n$* (this top-degree part has been indicated in (1.5) by the curly brace)?

With this application in mind, the conjectural formula from [5] (which was presented in a rather sketchy way in Section 2.5) gives the following concrete prediction for $\mathrm{Ch}_n^{\mathrm{top}}$.

**Conjecture 2.1.**
$$\mathrm{Ch}_n^{\mathrm{top}} = \frac{-1}{n!} \sum_{M,\prec} \gamma^{n+1-|\mathcal{V}|} \mathfrak{N}_M, \tag{2.2}$$

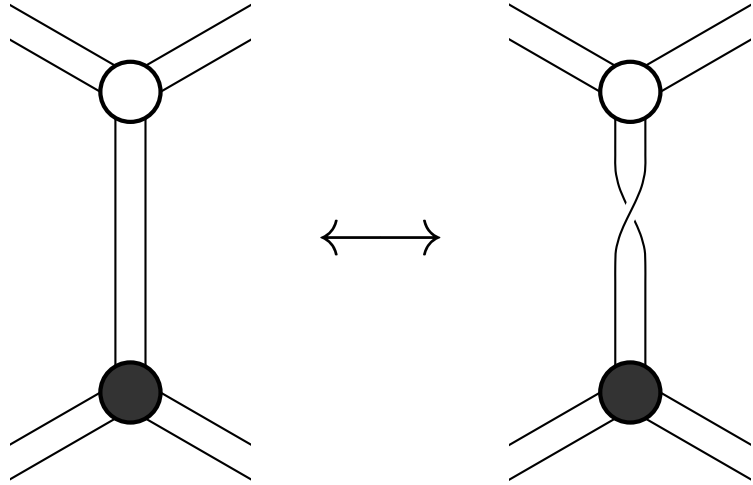

Figure 4: Twisting an edge of a ribbon graph.

*where the sum runs over **non-oriented** maps with n edges **and one face** with the following property: for each $1 \leq i < n$ if we remove the first i edges of the ribbon graph of M (according to the linear order $\prec$) then the resulting ribbon graph has the same number of connected components as the number of its faces.*

In the new context of $\mathrm{Ch}_n^{\mathrm{top}}$ we have to overcome two difficulties in order to fulfill our plan from Section 2.1: (A) we have to find an analogue of the abstract characterization of Jack characters $\mathrm{Ch}_n$ from Definition 1.2 which would work for $\mathrm{Ch}_n^{\mathrm{top}}$, and (B) we have reformulate the conjectural formula (2.2) into some more convenient form.

## 2.7 Untwisting the edges

The solution to the second difficulty (B) is provided by the following bijection.

**Theorem 2.2** ([3])**.** *There exists a bijection between:*

*(S1) the set of pairs $(M, \prec)$, where M is a **non-oriented** map with n edges **and one face**, and $\prec$ is a linear order on the set of its edges such that the condition from Conjecture 2.1 is fulfilled, and*

*(S2) the set of pairs $(M, \prec)$, where M is an **oriented** map with n edges **and arbitrary number of faces**, and $\prec$ is an **arbitrary** linear order on the set of its edges.*

*Sketch of a proof.* This bijection is performed by applying a twist (cf. Figure 4) to some selected edges of the ribbon graph of *M*. More specifially, we can view each ribbon graph which contributes to one of the above two sets as created from the void (i.e. a collection of isolated vertices) by adding ribbons in some specific order (the opposite of

≺). For each ribbon which is to be added we have to specify its 'orientation', i.e. one of the two ways in which it is to be attached; these two ways differ by the twist from Figure 4. On the other hand, we have to assure that the new ribbon is attached in a way that respects the condition related to a specific set (S1) or (S2). One can check that — regardless of the specific set (S1) or (S2) — for each position of the endpoints of a new ribbon there is either one choice (if the new ribbon is a bridge of a leaf) or there are two choices (otherwise) of a legitimate 'orientation'. In particular, the number of legitimate choices is the same, no matter if we construct a ribbon graph from (S1) or from (S2). This implies existence of a bijection. □

**Corollary 2.3.** *Conjecture 2.1 is equivalent to the following formula:*

$$\mathrm{Ch}_n^{\mathrm{top}} = (-1) \sum_M \gamma^{n+1-|\mathcal{V}|}\, \mathfrak{N}_M,$$

*where the sum runs over **oriented** maps with n edges and **arbitrary number of faces**.*

## 3 Towards the proof: abstract characterization of $\mathrm{Ch}_n^{\mathrm{top}}$

We need an abstract characterization of the top-degree part $\mathrm{Ch}_n^{\mathrm{top}}$ of Jack characters; a characterization which would use only intrinsic properties of $\mathrm{Ch}_n^{\mathrm{top}}$ and which would not refer to the much more complicated Jack character $\mathrm{Ch}_n$. Our characterization of $\mathrm{Ch}_n^{\mathrm{top}}$ is quite analogous to the characterization of $\mathrm{Ch}_n$ from (J1)–(J4). The only challenging task was to find a proper replacement for condition (J3) about the vanishing of the characters on small Young diagrams. Indeed, since the difference

$$\delta_n := \mathrm{Ch}_n - \mathrm{Ch}_n^{\mathrm{top}}$$

between the Jack character and its top-degree part is usually non-zero, if in (J3) we mechanically replace the Jack character $\mathrm{Ch}_n$ by its top-degree part $\mathrm{Ch}_n^{\mathrm{top}}$, we would get a statement which is clearly false.

A solution which we present in Lemma 3.3 is to require that certain *linear combinations* (over $|\lambda| < n$) of the values of $\mathrm{Ch}_n^{\mathrm{top}}(\lambda)$ vanish. These linear combinations were chosen in such a way that analogous linear combinations for $\delta_n$ vanish tautologically for any polynomial function $\delta_n$ which is of degree smaller than the degree of $\mathrm{Ch}_n$.

### 3.1 Top degree part $\mathrm{Ch}_n^{\mathrm{top}}$ revisited

The *filtration* on the algebra of polynomial functions which was defined Definition 1.1 may be refined to a certain natural *gradation* which corresponds to homogenous Stanley character polynomials of specified degrees. It makes sense therefore to define the top-degree $\mathrm{Ch}_n^{\mathrm{top}}$ as the homogeneous part of $\mathrm{Ch}_n$, of degree $d+1$, cf. (1.5).

In order to avoid introducing this gradation, the top-degree of $\mathrm{Ch}_n$ may defined as an *arbitrary* function $\mathrm{Ch}_n^{\mathrm{top}}\colon \mathbb{Y} \to \mathbb{Q}\left[A, A^{-1}\right]$ with the property that $\mathrm{Ch}_n - \mathrm{Ch}_n^{\mathrm{top}}$ is a polynomial function of degree at most $n$. With this definition, $\mathrm{Ch}_n^{\mathrm{top}}$ is defined only up to terms of degree $n$.

## 3.2  Operations on functions on $\mathbb{Z}^\ell$ and $\mathbb{Y}$

*Definition* 3.1. If $F = F(\lambda_1, \ldots, \lambda_\ell)$ is a function of $\ell$ arguments and $1 \leq j \leq \ell$, we define a new function $\Delta_{\lambda_j} F$ by

$$\left(\Delta_{\lambda_j} F\right)(\lambda_1, \ldots, \lambda_\ell) := F(\lambda_1, \ldots, \lambda_{j-1}, \lambda_j + 1, \lambda_{j+1}, \ldots, \lambda_\ell) - F(\lambda_1, \ldots, \lambda_\ell).$$

Any function $F$ on the set of Young diagrams can be viewed as a function $F(\lambda_1, \ldots, \lambda_\ell)$ defined for all non-negative integers $\lambda_1 \geq \ldots \geq \lambda_\ell$. We will extend its domain.

*Definition* 3.2. If $(\xi_1, \ldots, \xi_\ell)$ is a sequence of non-negative integers, we denote

$$F^{\mathrm{sym}}(\xi_1, \ldots, \xi_\ell) := F(\lambda_1, \ldots, \lambda_\ell),$$

where $(\lambda_1, \ldots, \lambda_\ell) \in \mathbb{Y}$ is the sequence $(\xi_1, \ldots, \xi_\ell)$ sorted in reverse order $\lambda_1 \geq \ldots \geq \lambda_\ell$.

## 3.3  Abstract characterization of $\mathrm{Ch}_n^{\mathrm{top}}$

The following result takes advantage of the deformation parameter $A$ on which Jack characters depend implicitly.

**Lemma 3.3.** *The top-degree $\mathrm{Ch}_n^{\mathrm{top}}$ can be characterized as a polynomial function $G\colon \mathbb{Y} \to \mathbb{Q}\left[A, A^{-1}\right]$ of degree at most $d := n+1$ for which an analogue of condition (J4) holds true and such that the equality*

$$[A^{n+1-2k}]\Delta_{\lambda_1} \cdots \Delta_{\lambda_k} G^{\mathrm{sym}}(\lambda_1, \ldots, \lambda_k) = 0 \tag{3.1}$$

*holds true for all $k \geq 0$ and all Young diagrams $\lambda \in \mathbb{Y}$ with at most $k$ rows and $|\lambda| \leq n - 1 - k$.*

*A substitute of a proof.* We start with a simple observation that $G := \mathrm{Ch}_n$ trivially fulfills the system of equations (3.1) because the left-hand side involves only the values of $G$ on the Young diagrams $\mu$ such that $|\mu| < n$ which vanish by the defining condition (J3).

Following Definition 1.1, $G$ is of the form

$$G(\lambda) = \underbrace{p_0(\gamma)}_{G_0(\lambda):=} + \underbrace{\sum_{\square_1 \in \lambda} p_1(\gamma, c_1)}_{G_1(\lambda):=} + \underbrace{\sum_{\square_1, \square_2 \in \lambda} p_2(\gamma, c_1, c_2)}_{G_2(\lambda):=} + \cdots. \tag{3.2}$$

The remaining difficulty is to show that:

(G1) the system of equations (3.1) *does not* involve the subdominant terms of the polynomials $p_0, p_1, \dots$; or, in other words, it is indeed fulfilled by

$$G := \text{Ch}_n^{\text{top}} = \text{Ch}_n + (\text{arbitary polynomial function of degree at most } n);$$

(G2) (3.1) is a system of equations which (together with condition (J4)) determines uniquely the top-degree parts of the polynomials $p_0, p_1, \dots$; or, in other words, it determines the top-degree part of $\text{Ch}_n^{\text{top}}$ uniquely.

The initial degree bound (J1) from Definition 1.1 implies that $[A^{n+1-2k}]G$ depends only on the polynomials $p_0, \dots, p_k$; in other words (3.1) can be seen as a kind of an upper-triangular system of equations over $k = 0, 1, 2, \dots$.

On the other hand, the finite difference calculus *seems to imply* that the contribution of the polynomials $p_0, \dots, p_{k-1}$ to $\Delta_{\lambda_1} \cdots \Delta_{\lambda_k} G^{\text{sym}}$ tautologically vanishes (*"the iterated finite difference $\Delta_{\lambda_1} \cdots \Delta_{\lambda_k}$ of a function of separated variables which does not depend on some variable is zero"*); in particular *it seems that* the left-hand side of (3.1) depends only on $p_k^{\text{top}}$, the homogeneous part of the polynomial $p_k$ of degree $n + 1 - 2k$:

$$[A^{n+1}]G^{\text{sym}}(\emptyset) = [A^{n+1}]p_0(\gamma) = p_0^{\text{top}}(-1), \tag{3.3}$$

$$[A^{n-1}]\Delta_{\lambda_1} G^{\text{sym}}(\lambda_1) = [A^{n-1}]p_1^{\text{top}}(-1, \lambda_1 + 1), \tag{3.4}$$

$$[A^{n-3}]\Delta_{\lambda_1}\Delta_{\lambda_2} G^{\text{sym}}(\lambda_1, \lambda_2) = 2[A^{n-3}]p_2^{\text{top}}(-1, \lambda_1 + 1, \lambda_2 + 1), \tag{3.5}$$

$$\vdots$$

It follows, in particular, that we have reached our goal (G1).

Consider the first nontrivial case $k = 1$ and the corresponding polynomial

$$c_1 \mapsto p_1^{\text{top}}(-1, c_1) \tag{3.6}$$

in a single variable; a polynomial which is of degree $n - 1$. Equations (3.1) and (3.4) provide the values of this polynomial in $n - 1$ points; furthermore condition (J4) specifies the top-degree coefficient of this polynomial. The above data determines (3.6) uniquely; it follows that the homogeneous polynomial $p_1^{\text{top}}$ is also uniquely determined which is a step towards our goal (G2).

Regretfully, the above considerations **are not correct for** $k \geq 2$. The source of the difficulties lies in the fact that (3.5) **does not hold true if** $\lambda_1 = \lambda_2$. Indeed, in this special case we encounter the subtlety related to the way the symmetrization works:

$$\Delta_{\lambda_1}\Delta_{\lambda_2} G^{\text{sym}}(\lambda_1, \lambda_2)|_{\lambda_1 = \lambda_2} =$$
$$G^{\text{sym}}(\lambda_1 + 1, \lambda_1 + 1) - G^{\text{sym}}(\lambda_1 + 1, \lambda_1) - G^{\text{sym}}(\lambda_1, \lambda_1 + 1) + G^{\text{sym}}(\lambda_1, \lambda_1) =$$
$$G(\lambda_1 + 1, \lambda_1 + 1) - 2G(\lambda_1 + 1, \lambda_1) + G(\lambda_1, \lambda_1).$$

It follows that the argument with the finite difference calculus does not work and, unexpectely, $\Delta_{\lambda_1}\Delta_{\lambda_2} G_1^{\text{sym}} \not\equiv 0$, where $G_1$ is defined by (3.2). The correct version of (3.5) for $\lambda_1 = \lambda_2$ turns out to be:

$$[A^{n-3}]\Delta_{\lambda_1}\Delta_{\lambda_2} G^{\text{sym}}(\lambda_1, \lambda_2)\Big|_{\lambda_1=\lambda_2} = 2[A^{n-3}]p_2^{\text{top}}(-1, \lambda_1+1, \lambda_2+1) - \frac{\partial}{\partial \lambda_1} p_1^{\text{top}}(-1, \lambda_1+1).$$

This kind of argument can be carried out for arbitrary $k \geq 1$; it turns out that the left-hand side of (3.1) involves only the values of the polynomials $p_l^{\text{top}}(1, \cdot)$ over $l \leq k$. The values of $p_k^{\text{top}}(1, \cdot)$ turn out to be specified in sufficiently many points to determine the polynomial $p_k^{\text{top}}$ recursively, in terms of $p_1^{\text{top}}, \ldots, p_{k-1}^{\text{top}}$. In this way one can show that we have reached our goals (G1) and (G2). □

## 4 The main result: closed formula for $\text{Ch}_n^{\text{top}}$

It is time for the second easy step (2) from Section 2.1: to verify that the formula from Corollary 2.3 indeed fulfills the characterization of $\text{Ch}_n^{\text{top}}$ given by Lemma 3.3.

**Theorem 4.1** ([11]).

$$\text{Ch}_n^{\text{top}} = (-1)\sum_M \gamma^{n+1-|\mathcal{V}|}\,\mathfrak{N}_M, \tag{4.1}$$

*where the sum runs over **oriented** maps with $n$ edges and **arbitrary number of faces**.*

*Sketch of a proof.* In order to prove (J2) we one can use some general methods for verifying that a given linear combination of bipartite graphs defines a polynomial function [6] in the sense of Definition 1.1.

The limit bounds (J1) on the degrees of the content polynomials as well as the normalization (J4) follow very easily from the asymptototics of the right-hand side of (4.1) in the limit as $\lambda \to \infty$ and $\gamma \to \infty$.

In order to prove that condition (3.1) holds true we show that the contribution of non-injective embeddings (which map two edges of the map into the same box of $\lambda$) vanishes. This can be done by finding a sign-reversing involution on the set of maps which are compatible with a prescribed embedding of the edges. Since there are no injective embeddings of a map into a Young diagram with a small number of boxes, it follows that (3.1) indeed holds true. □

## 5 Outlook: closed formula for sub-topdegree part of $\text{Ch}_n$?

Can we repeat our success story and find a closed formula for, say, the sub-topdegree part of Jack character $\text{Ch}_n$? Unfortunately, Theorem 4.1 does not offer any hints how such a formula could look like. On the bright side, the incorrect formula from [5] might give some hints. The future will show.

# Acknowledgements

Research supported by *Narodowe Centrum Nauki*, grant number 2014/15/B/ST1/00064.